\documentstyle[amscd]{amsart}
\newtheorem{thm}{Theorem}[section]
\newtheorem{lm}[thm]{Lemma}
\newtheorem{cor}[thm]{Corollary}
\newtheorem{pro}[thm]{Proposition}

\newtheorem{th}{Theorem}

\theoremstyle{definition}
\newtheorem{df}[thm]{Definition}

\theoremstyle{remark}

\numberwithin{equation}{section}

\def \R {\bold{R}}
\def \Z {\bold{Z}}

\def \CR {\cal R}
\def \CS {\cal S}
\def \J {\cal J}

\def \CM {\cal M}

\def \re {representation}
\def \hd {Heegaard decomposition}

\def \a {\alpha}
\def \b {\beta}
\def \c {\chi}
\def \g {\gamma}
\def \G {\Gamma}

\def \e {\eta}

\def \Lam {\Lambda}
\def \lam {\lambda}
\def \lc {\lambda_{CL}}

\def \n {\nabla}
\def \o {\omega}
\def \p {\phi}
\def \r {\rho}
\def \s {\sigma}

\def \th {\theta}
\def \ti {\tilde}

\def \O {\Omega}

\def \bd {\partial}

\def \pjh {\overline{\partial}_{J,H}}

\def \x {\times}
\def \ox {\otimes}
\def \ve {\varepsilon}
\def \D {\Delta}

\def \CP {{\bold{CP}}^1}

\begin{document}

\baselineskip.525cm

\title[Casson-Lin's invariant of a knot, Floer homology]
{Casson-Lin's invariant of a knot and Floer homology}

\author[Weiping Li]{Weiping Li}
\address{Department of Mathematics, Oklahoma State University \newline
\hspace*{.375in}Stillwater, Oklahoma 74078-0613}
\email{wli@@math.okstate.edu}
\thanks{Partially supported by a summer research award from the
College of Arts and Sciences at Oklahoma State University.}

\begin{abstract}
A. Casson defined an intersection number 
invariant which can be roughly thought of as the number of conjugacy classes of
irreducible representations of $\pi_1(Y)$ into $SU(2)$ counted with signs,
where $Y$ is an oriented integral homology 3-sphere.
X.S. Lin defined an similar invariant (signature of a knot)
to a braid representative of a knot in
$S^3$. In this paper, we give a
natural generalization of the Casson-Lin's invariant to be 
(instead of using the instanton Floer homology) 
the symplectic Floer homology for
the representation space (one singular point) of $\pi_1(S^3 \setminus K)$
into $SU(2)$ with trace-free along all meridians. The symplectic
Floer homology of braids is a new invariant of knots 
and its Euler number of such
a symplectic Floer homology is the negative of the Casson-Lin's invariant.
\end{abstract}

\maketitle

\section{Introduction}

Around 1986 there appeared two invariants for a closed, oriented, 3-manifold
$Y$ with the integral homology of $S^3$.
The first, introduced by Casson (see \cite{wa}
), is an integer-valued invariant which roughly counts the number of irreducible
$SU(2)$-representations of the fundamental group $\pi_1(Y)$. The second is a
$Z_8$-graded homology theory developed by Floer in \cite{fl}, and is based upon
an application of the Morse theory to the Chern-Simons functional on
the space ${\cal B}_Y$ of equivalent
classes of $SU(2)$ connections on $Y$. 
The two invariants are closely related because
irreducible $SU(2)$-representations can be interpreted as the critical points
of the Chern-Simons functional. The Euler characteristic of the instanton
Floer homology is twice Casson's invariant of the integral homology
3-sphere $Y$.

Analogous to Casson invariant, Lin \cite{lin} constructed an invariant of
representation spaces corresponding to a braid representative of a knot $K$
in $S^3$. It turns out that the Lin's invariant of a knot is the 
signature of the knot. In \cite{lin}, Lin used the special {\re}s of
$\pi_1( S^3 \setminus K)$ into $SU(2)$ such that all meridians of $K$ are
represented by trace zero matrices in $SU(2)$. It naturally arises two questions
to understand more on the Casson-Lin's invariant of a knot. (1) Does there
exists a similar Floer homology generalization for the
Casson-Lin's invariant of a knot ? (2) What kind of invariant one can get
by replacing the trace zero to a fixed trace {\re}s of knot groups ?
Recently, Cappell, Lee and Miller have studied the second question 
from the symplectic theory point of view in \cite{clm}
and defined a so-called equivariant Casson invariant for 3-manifold with
cyclic group action. Independently, Herald studied the same problem from
the gauge theory point of view in \cite{he}.
The present paper is to give a Floer homology type of generalization for
the Casson-Lin's invariant of a knot for the first question.

In this paper, we derive the symplectic Floer homology, rather than
the instanton Floer homology, as the generalization of the
Casson-Lin's invariant of a knot. It is not clear how to make the
gauge theory point of view to realize that
the critical points of certain functional (Chern-Simons functional)
will correspond to the {\re} of $S^3 \setminus K$ with trace free 
along all meridians.
Instead of the instanton Floer homology, the symplectic Floer homology
also provides the generalization of the Casson invariant which help us to 
shift the point of view. From the Atiyah conjecture as described
in \cite{ll}, \cite{ll2}, we have seen that
there is really no difference between the instanton Floer homology
and the symplectic Floer homology for integral homology 3-spheres.
So we can construct a generalization for the
Casson-Lin's invariant of a knot by constructing the symplectic Floer homology.

The paper is organized as follows. In \S 2, we recall the 
Casson-Lin's invariant of a knot from the symplectic topology point of
view which also provides some basic ingredients to construct the symplectic
Floer homology. In \S 3, we mainly concentrate on the reducible {\re}
of $\pi_1(S^3 \setminus K)$ and discuss the Walker-type correction term
for this isolated reducible $U(1)$-{\re}. We briefly review the symplectic
Floer homology in \S 4.1 and show that the symplectic Floer homology is an 
invariant of a knot in \S 4.2. Some simple calculations of the new invariant
are given in \S 5 in order to fix the sign between the Euler number and the
invariant $\lc$.

\section{Casson-Lin's invariant of a knot}

In this section, we recall Lin's construction in \cite{lin}
for the
intersection number of the representation spaces corresponding to a braid
representative of a knot $K$ in $S^3$. 
The presentation of the Casson-Lin's invariant here 
is from the symplectic
geometry point of view in order to extend the Casson-Lin's invariant to the
symplectic Floer homology in \S 4.

Let $(S^3, D^3_+, D^3_-, S^2)$ be a Heegaard decomposition of $S^3$ with genus
$0$, where
\[ S^3 = D^3_+ \cup_{S^2} D^3_- , \ \ \ \bd D^3_+ = \bd D^3_- =
D^3_+ \cap D^3_- = S^2 .\]
Suppose that a knot $K \subset S^3$ is 
in general position with respect to this Heegaard decomposition. So
$K \cap S^2 = \{x_1, \cdots, x_n, y_1, \cdots, y_n\}$, $K \cap D^3_{\pm }$
is a collection of unknotted, unlinked arcs
$\{ \g_1^{\pm }, \cdots, \g_n^{\pm } \} \subset D^3_{\pm }$, 
where $\bd \g_i^- = \{x_i, y_i\}$ and $\{\g_1^+, \cdots, \g_n^+\} = K \cap
D^3_+$ becomes a braid of $n$ strands inside $D^3_+$. Denote a word $\b $ in the
braid group $B_n$. For the top end point $x_i$ of $\g_i^+$, the bottom 
end points of $\{\g_1^+, \cdots, \g_n^+\}$ gives a permutation of 
$\{y_1, \cdots , y_n\}$ which generates a map 
\[ \pi : B_n \to S_n ,\] 
where $\pi (\b )$ is the permutation of $\{y_1, \cdots , y_n\}$ in the symmetric
group of $n$ letters. Let $K = \overline{\b }$ be the closure of $\b$.
It is well-known that there is a correspondence between a knot and a braid $\b$
with $\pi (\b )$ is a complete cycle of the $n$ letters (see \cite{bi}).

There is a corresponding {\hd} for the complement of a knot $K$,
\[ S^3 \setminus K = (D^3_+ \setminus K) \cup_{(S^2 \setminus K)} (D^3_-
\setminus K ) , \]  
\[  D^3_{\pm } \setminus K = D^3_{\pm } \setminus (D^3_{\pm }
\cap K), \ \ \ S^2 \setminus K = S^2 \setminus (S^2 \cap K) . \]
Thus from Seifert-Van Kampen theorem,
\[ \begin{array} {ccc}
\pi_1 (S^2 \setminus K) & \to & \pi_1 (D^3_+ \setminus K)\\
\downarrow & & \downarrow \\
\pi_1 (D^3_- \setminus K ) & \to & \pi_1 (S^3 \setminus K) , 
\end{array} \]
we obtain a pull-back of {\re} spaces
\begin{equation} \label{pullback}
 \begin{array} {ccc}
{\CR}(S^2 \setminus K)& \leftarrow & {\CR}(D^3_+ \setminus K)\\
\uparrow & & \uparrow \\
{\CR}(D^3_- \setminus K ) & \leftarrow & {\CR} (S^3 \setminus K) , 
\end{array} \end{equation}
where ${\CR }(X) = Hom (\pi_1 (X), SU(2))/SU(2)$ for $X = 
S^2 \setminus K, D^3_{\pm } \setminus K, S^3 \setminus K$.
The $SU(2)$ admits a natural Riemannian metric which arises by translation
from the inner product on the Lie algebra of $S^3$. The natural identification
of $S^3$ with $SU(2)$ is an isometry of the standard metric on $S^3$ and the
above Riemannian metric on $SU(2)$. Also note that $S^3 \setminus K$ is an
Eilenberg-MacLane space $K(\pi_1(S^3 \setminus K), 1)$ for
classical knot $K$ due to the asphericity theorem of Papakyriakopoulos.
\medskip

In \cite{ma}, Magnus used the trace free matrices to represent the
generators of a free group to show that the faithfulness of a
{\re} of braid groups in the automorphism groups of the rings generated
by the character functions on free groups. This is original idea 
to have {\re}s with trace free along all meridians which Lin worked in
\cite{lin} to define the knot invariant. It has been carried out by
Cappell, Lee and Miller in \cite{clm} for the {\re} of knot groups with
the trace of the meridian fixed (not necessary zero). We only give the
generation to symplectic Floer homology of the Casson-Lin's invariant for trace
free {\re}s of knot groups. For general case as in \cite{clm}, we will discuss
the corresponding Floer homology elsewhere.
\medskip

Let $(M, \o)$ be a $2n$-dimensional symplectic manifold and
$\p : M \to M$ be a symplectic diffeomorphism, i.e. $\o$ is a nondegenerate
closed 2-form and $\p^* \o = \o$. By choosing an almost complex structure $J$ on
$(M, \o)$ such that $\o (\cdot, J \cdot)$ defines a Riemannian metric,
we have an integer valued second cohomology class $c_1(M) \in H^2(M, Z)$
(the first Chern class). Note that the space of all almost complex structures
which are compatible with $\o $ is connected, so $c_1(M)$ is uniquely determined
by $\o$. There are two homomorphisms
\[I_{\o }: \pi_2(M) \to {\R }; \ \ \ \ \ \
I_{c_1}: \pi_2(M) \to {\Z }. \]
\begin{df} \label{monotone}
The symplectic manifold $(M, \o)$ is called monotone if there exists
a nonnagetive $\a \geq 0$ such that
$I_{\o } = \a I_{c_1}$ on $\pi_2(M)$ or $\pi_2(M) = 0$.
See also \cite{fl3}.
\end{df}

Let ${\CR}(S^3 \setminus K)^{[i]}$ be the 
space of $SU(2)$ {\re}s $\r : \pi_1 (S^3 \setminus K) \to SU(2)$
such that
\begin{equation} \label{trace}
 \r ([m_{x_i}]) \sim \left( \begin{array}{cc}
i & 0 \\ 0 & -i \end{array} \right) , \ \ \ 
\r ([m_{y_i}]) \sim \left( \begin{array}{cc}
i & 0 \\ 0 & -i \end{array} \right) , \end{equation}
where $m_{x_i}, m_{y_i}, i =1, 2, \cdots, n$ are the meridian circles around
$x_i, y_i$ respectively. Note that $\pi_1(S^3 \setminus K)$ is generated by
$m_{x_i}, m_{y_i}, i =1, 2, \cdots, n$
and one relation $\prod^n_{i=1} m_{x_i} = \prod^n_{i=1} m_{y_i}$.
Corresponding to (\ref{pullback}), we have
\begin{equation} \label{irre}
 \begin{array} {ccc}
{\CR}(S^2 \setminus K)^{[i]} & \leftarrow & {\CR}(D^3_+ \setminus K)^{[i]}\\
\uparrow & & \uparrow \\
{\CR}(D^3_- \setminus K )^{[i]} & \leftarrow & {\CR} (S^3 \setminus K)^{[i]} .
\end{array} 
\end{equation}
The conjugacy class in $SU(2)$ is completely determined by its trace.
So the condition (\ref{trace}) can be reformulated for $\r \in {\CR}(X)^{[i]}$,
\begin{equation} \label{22}
trace \r ([m_{x_i}]) = trace \r ([m_{y_i}]) = 0 . \end{equation}
The space ${\CR}(S^2 \setminus K)^{[i]}$ can be identified with the space
of $2n$ matrices $X_1 \cdots, X_n, Y_1, \cdots, Y_n$ in $SU(2)$ satisfying
\begin{equation} \label{zero}
trace (X_i) = trace (Y_i) = 0, \ \ \ \ \mbox{for $i =1, \cdots, n$},
\end{equation}
\begin{equation} \label{product}
X_1 \cdot X_2 \cdots X_n = Y_1 \cdot Y_2 \cdots Y_n .
\end{equation}

\begin{lm} \label{qmono}
Let $Q_n$ be the space $\{(X_1, \cdots , X_n) \in SU(2)^n | \ \ 
trace (X_i) = 0, i =1, \cdots, n \}$.
Thus $Q_n$ is a monotone symplectic manifold of dimension $2n$.
\end{lm}
Proof: An element in $SU(2)$ can be viewed as $\left[ \begin{array}{cc}
a+bi & c+di \\ -c +d i & a - bi \end{array} \right]$, with $a^2 + b^2 + c^2
+ d^2 = 1$. A matrix $X_j$ in $Q_n$ is a matrix with $trace (X_j) = 2 a_j = 0$.
I.e. the set of all such matrices equal to  $\{(b_j, c_j, d_j) | 
b_j^2 +c_j^2 + d_j^2 = 1 \}$. Hence $Q_n$ is the product
$S^2 \x S^2 \cdots \x S^2$ of $n$ copies of 2-sphere with radius $1$.
So all $S^2$'s have the same area $4 \pi$. It is known that
$S^2 \x S^2 \cdots \x S^2$ is a monotone symplectic manifold if and only if
all $S^2$ has the same area (see \cite{fl3} page 577).  
Thus the result follows.  \qed

It is well known that there is a natural symplectic structure on the coadjoint
orbits by Kostant-Souriau theorem. If $H^1({\bf g}) = H^2({\bf g}) = 0$
for Lie algebra ${\bf g}$ of a Lie group $G$, every $\o \in Z^2({\bf g})$
is exact, $\o = d \a$. The $\a$ is unique due to $H^1({\bf g}) =0$. There
is one-to-one correspondence between $G$-orbits 
$orb_G(\o) = \{ad^*_a \o : a \in G\}$ in $Z^2({\bf g})$ and
$G$-orbits $orb_G(\a ) = \{ad^*_a \a : a \in G\}$ in ${\bf g}^*$, where
$ad^*$ is the coadjoint action and ${\bf g}^*$ is the dual space of ${\bf g}$.
The symplectic form $\overline{\o}$ on $orb_G(\a ) \cong G/H_{\o}$ 
can be expressed by
\[ \overline{\o} (X_{\a}, Y_{\a}) = - \a ([X, Y]), \]
for $\pi : G \to G/H_{\o}$ and $X, Y \in {\bf g}, X_{\a} = \pi_*X, 
Y_{\a} = \pi_* Y$.
For $G = SU(2)$, 
\[ orb_G(B) = \{ad^*_a B : a\in SU(2)\} \ \ \ \ \mbox{for $B= \left( 
\begin{array}{cc} 
i &  \\
 & - i \end{array} \right).$ } \]
is an orbit $\{i a^{-1}Ba : a \in SU(2) \}$ in ${\bf su}(2)^*$. Note that 
$a^{-1}Ba = B$ if and only if $a = \left( \begin{array}{cc}
t & \\ 
 & t^{-1} \end{array} \right) $ with $t \in U(1)$. So
\[ orb_G(B) = SU(2)/U(1) \cong {\CP} .\]
A symplectic form on ${\CP}$ is given by  $\o (X, Y) = (- B, [A^X, A^Y])$ for
$\pi_*A^X = X, \pi_* A^Y = Y$. The trace free condition ensures to lie
in the orbit $orb_G(B)$. Hence $Q_n$ can be naturally identified
with ${\CP} \x \cdots \x {\CP}$ as a product of $n$ copies of
${\CP}$, again a monotone symplectic manifold.
\medskip

Denote $\ti{S}^2$ be the branched double covering of $S^2$ with $K \cap S^2$
as branched set. Then $Z_2$ operates on ${\CR }(\ti{S}^2)$ and the 
component of the fixed point set ${\CR }^*(\ti{S}^2)$ (the
irreducible $SU(2)$ {\re}s) with traceless condition coincides with
${\CR}^*(S^2 \setminus K)^{[i]}$ the irreducible $SU(2)$ {\re}s with
traceless condition. By \cite{go}, the irreducible components carry the
natural symplectic structure and dimension $4n-6$ has been verified
in \cite{lin} p342- 344.

\begin{lm} \label{mos}
The space ${\CR}^*(S^2 \setminus K)^{[i]}$ is a monotone symplectic
manifold of dimension $4n-6$.
\end{lm}
Proof: Note that ${\CR}^*(S^2 \setminus K)^{[i]} = (H_n \setminus S_n)/SU(2)$
in Lin's notation \cite{lin}, where
\[H_n = \{(X_1, \cdots, X_n, Y_1, \cdots, Y_n) \in Q_n \x Q_n | \ \
X_1 \cdots X_n = Y_1 \cdots Y_n \},\]
$S_n$ is the subspace of $H_n$ consisting of all the reducible points. 
Here $H_n \setminus S_n$ is the total space of a $SU(2)$-fiber bundle over
${\CR}^*(S^2 \setminus K)^{[i]}$. It is remarkable that the points in 
$H_n$ where $SU(2)$-action is not locally free are precisely the set $S_n$.
Note that any product of monotone
symplectic manifolds is again monotone. So
we have $I_{\o} = \a I_{c_1}$ for $\a \geq 0$ over $(Q_n \x Q_n, \o)$.
$(H_n, \o)$ is pre-symplectic space which is degenerated along
the $SU(2)$-fiber.  
For $n \geq 2$, $S_n$
has codimension bigger that $2$.
The tangent space of $H_n \setminus S_n$ splits into a direct sum of
tangent spaces along vertical and horizontal directions,
\begin{equation} \label{split}
T(H_n \setminus S_n) = T {\CR}^*(S^2 \setminus K)^{[i]} \oplus TSU(2).
\end{equation} 
Here $TSU(2) \cong T^*SU(2)$ carries a canonical symplectic form
$\o_0 = - d \lam $, where $\lam$ is the Liouville 1-form.
Thus $\o - \o_0$ is nondengenerate closed 2-form on
${\CR}^*(S^2 \setminus K)^{[i]}$. Let $\pi : H_n \setminus S_n \to 
{\CR}^*(S^2 \setminus K)^{[i]}$ be the $SU(2)$-orbit projection, so 
$\o - \o_0 = \pi_*(\o)$. From the fibration and (\ref{split}), 
the Levi-Civita connection $A$ on $T(H_n \setminus S_n)$ is essentially
the Levi-Civita connection $a$ on $T {\CR}^*(S^2 \setminus K)^{[i]}$
since $SU(2)$ is parallelizable, $A = a + \th$, where $\th$ 
is the trivial connection.
Thus by Chern-Weil theorem, $2 \pi i c_1 = F_A$, we have
$\pi^* c_1 = c_1 |_{H_n \setminus S_n} \in H^2(H_n \setminus S_n; Z)$, i.e.
$c_1|_{H_n \setminus S_n}$ is the natural pullback of the 
first Chern class $c_1$ in
$H^2({\CR}^*(S^2 \setminus K)^{[i]}, {\Z})$. 
By the homotopy exact sequence of the fibration $\pi: H_n \setminus S_n
\to {\CR}^*(S^2 \setminus K)^{[i]}$ with fiber $SU(2)$,
\[0 = \pi_2 (SU(2)) \to \pi_2(H_n \setminus S_n) \stackrel{\pi_*}{\rightarrow}
\pi_2({\CR}^*(S^2 \setminus K)^{[i]}) \to \pi_1(SU(2)) = 0 , \]
we obtain the natural isomorphism induced by $\pi_*$
$\pi_2(H_n \setminus S_n) \cong 
\pi_2({\CR}^*(S^2 \setminus K)^{[i]})$.  
By Stokes theorem, 
$I_{\o_0} = 0$. Hence for any 
$v \in \pi_2 ({\CR}^*(S^2 \setminus K)^{[i]}) = \pi_2(H_n \setminus S_n)$ (after
identification by $\pi_*$), 
\[I_{\o}(v) = I_{\pi_*(\o)}(v) + I_{\o_0}(v) = I_{\pi_*(\o)}(v) . \]
\[ I_{c_1}(v)  = I_{\pi^*(c_1)}(v).\]
Thus $I_{\pi_* \o} = \a I_{c_1}$ for $\a \geq 0$ on
${\CR}^*(S^2 \setminus K)^{[i]}$. Thus we get the monotonicity
for ${\CR}^*(S^2 \setminus K)^{[i]}$. \qed

Given $\b \in B_n$, we denote by $\G_{\b}$ the graph of $\b $ in $Q_n \x Q_n$,
i.e.
\[\G_{\b} = \{(X_1, \cdots, X_n, \b (X_1) , \cdots, \b (X_n)) \in Q_n \x Q_n\}.
\]
As an automorphism of the free group $Z[m_1]*Z[m_2]* \cdots *Z[m_n]$, this
element $\b \in B_n$ preserves the word $[m_1] \cdots [m_n]$. Thus we have
\[ X_1 \cdots X_n = \b (X_1) \cdots \b (X_n), \]
or in other words $\G_{\b}$ is a subspace of $H_n$. In fact, for $\overline{\b}
= K$, this subspace $\G_{\b}$ coincides with  the subspace of {\re}s
$\r : \pi_1(S^2 \setminus K) \to SU(2)$ in $H_n$ which can be extended
to $\pi_1(D^3_+ \setminus K)$, 
$\G_{\b} = Hom (\pi_1(D^3_+ \setminus K), SU(2))^{[i]}.$ 
Hence the space ${\CR}^*(D^3_+ \setminus K)^{[i]} = \G_{\b, irre}/SU(2)$ 
is the irreducible $SU(2)$ {\re}s with traceless condition over
$D^3_+ \setminus K$.

In the special case $\b = id$, then $\G_{id}$ represents the diagonal in
$Q_n \x Q_n$,
\[\G_{id} = \{ (X_1, \cdots, X_n, X_1, \cdots, X_n) \in Q_n \x Q_n\} . \]
Since $K \cap D^3_-$ represents the trivial braid, this space $\G_{id}
\subset H_n$ can be identified with the subspace of {\re}s
in $Hom(\pi_1(S^2 \setminus K), SU(2))^{[i]}$ can be extended to
$\pi_1(D^3_- \setminus K)$, i.e.
\[ \G_{id} = Hom (\pi_1(D^3_- \setminus K), SU(2))^{[i]} . \]
By Seifert, Van-Kampen Theorem, the intersection
$\G_{\b} \cap \G_{id}$ is the same as the space of {\re}s of
$\pi_1(S^3 \setminus K)$ satisfying the monodromy condition $[i]$ 
(see (\ref{pullback})),  
\[\G_{\b} \cap \G_{id} = Hom (\pi_1(S^3 \setminus K), SU(2))^{[i]} . \]

Given $\b \in B_n$ with $\overline{\b } = K$, there is an induced diffeomorphism
(still denoted by $\b $) from $Q_n$ to itself. Such a diffeomorphism also
induces a diffeomorphism $\p_{\b }:  {\CR}^*(S^2 \setminus K)^{[i]} \to 
{\CR}^*(S^2 \setminus K)^{[i]}$.
\begin{lm} \label{symdif}
For $\b \in B_n$ with $\overline{\b } = K$, the induced diffeomorphism
$\p_{\b }:  {\CR}^*(S^2 \setminus K)^{[i]} \to 
{\CR}^*(S^2 \setminus K)^{[i]}$ is symplectic, and the fixed point
set of $\p_{\b }$ is ${\CR}^*(S^3 \setminus K)^{[i]}$.
\end{lm}
Proof: For $\b \in B_n$ with $\overline{\b } = K$, there is an induced
diffeomorphism (see \cite{lin}). Note that $\b $ maps each $X_j$ to a conjugate
of some $X_{j^{'}}$, it leaves $Q_n$ invariant, so it gives
rise an area-preserving diffeomorphism of $Q_n$. I.e. 
\[ Area(S^2_j) = Area(S^2_{j^{'}}) = Area(S^2_{\b (j)}) ,\]
where $S^2_j = \{ X_j \in SU(2) | tr X_j = 0\}$.
Let $\o = \sum_{j=1}^n \o_j$ be the symplectic form on $Q_n$ with
the symplectic form $\o_j$ on the $j$-th 2-sphere by Lemma~\ref{qmono}. 
Since $\b \in B_n$ is a complete cycle of $n$-letters, so
\[\b^*(\o ) = \sum _{j=1}^n \o_{\b (j)} = \o  .\]
From the identification
between ${\CR}^*(S^2 \setminus K)^{[i]}$ and $(H_n \setminus S_n)/SU(2)$
in Lemma~\ref{mos}, 
$\b^* \o_0 = \o_0$ on the $SU(2)$-orbit, the braid
$\b $ also induces a diffeomorphism $\p_{\b }$
which preserves the induced symplectic
structure on $(H_n \setminus S_n)/SU(2)$:
\begin{eqnarray*}
 \p_{\b }^*(\pi_* \o) & = & \b^* (\o - \o_0) \\
 &= & \b^* (\o) - \b^* (\o_0) \\
 & = & \o - \o_0  = \pi_* \o .
\end{eqnarray*}
Note that 
$\overline{\G}_{\b} = ({\G}_{\b} \setminus ({\G}_{\b} \cap S_n))/SU(2)$
is the graph of $\p_{\b }$. 
Similarly, $\overline{\G}_{id}$ can be thought as a ``diagonal''.
By Seifert, Van-Kempen theorem (\ref{irre}), it is clear that
\[ Fix (\p_{\b }|_{{\CR}^*(S^2 \setminus K)^{[i]}}) =
\overline{\G}_{\b} \cap \overline{\G}_{id} =
{\CR}^*(S^3 \setminus K)^{[i]} .\]
Thus we obtain the desired result. \qed

Let $\ti{D^3}_{\pm }$ be the double branched covering of ${D^3}_{\pm }$ with
${D^3}_{\pm } \cap K$ as branched set. There is an induced
$Z_2$-action on the {\re} space ${\CR }(\ti{D^3}_{\pm })$, and
$\overline{\G }_{\b} = {\CR}^*(D^3_+ \setminus K)^{[i]}$, 
$\overline{\G}_{id} = {\CR}^*(D^3_-\setminus K)^{[i]}$ 
can be regarded as components of
the fixed point set ${\CR }(\ti{D^3}_{\pm })^{Z_2}$. 
It follows that
${\CR}^*(D^3_+ \setminus K)^{[i]}, {\CR}^*(D^3_-\setminus K)^{[i]}$ 
have orientations
inherited from the orientation
on ${\CR }(D^3_{\pm })^{Z_2}$.
These oriented submanifolds 
${\CR}^*(D^3_+ \setminus K)^{[i]}, {\CR}^*(D^3_-\setminus K)^{[i]}$ 
intersects each other in a compact subspace of
${\CR}^*(S^2 \setminus K)^{[i]}$ from Lemma 1.6 in \cite{lin}. 
Hence we can perturb $\p_{\b }$ (via
Hamiltonian vector field if necessary), or
in another words perturb ${\CR}^*(D^3_+ \setminus K)^{[i]}$ to 
$\hat{{\CR}}^*(D^3_+ \setminus K)^{[i]}$ by a compactly support isotopy so that
$\hat{{\CR}}^*(D^3_+ \setminus K)^{[i]}$ intersects
${\CR}^*(D^3_-\setminus K)^{[i]}$ transversally at a finite number of 
intersection points. Denote the perturbed symplectic diffeomorphism
by $\ti{\p }_{\b }$. So its fixed points are all nondegenerated.

\begin{df} \label{cassonlin}
The Casson-Lin invariant of a knot $K = \overline{\b}$ is given by
counting the algebraic intersection number of 
$\hat{{\CR}}^*(D^3_+ \setminus K)^{[i]}$ and 
${\CR}^*(D^3_-\setminus K)^{[i]}$, or the algebraic number of
$Fix(\ti{\p }_{\b })$,
\[\lam_{CL} (K) =  \lam_{CL}(\b ) =  \# Fix(\ti{\p }_{\b }) = 
 \# (
\hat{{\CR}}^*(D^3_+ \setminus K)^{[i]} \cap {\CR}^*(D^3_-\setminus K)^{[i]}) .\]
\end{df}

The results proved in \cite{lin} show that the Casson-Lin invariant 
$\lam_{CL} (K) = \lam_{CL}(\b )$
is independent of its braid representatives, i.e. $\lam_{CL}(\b )$ is invariant
under the Markov moves of type I and type II on $\b$.

\section{Lagrangian perturbations on the singular {\re}}

In this section, we study the Walker type correction term around 
an isolated $U(1)$ reducible {\re} of the knot group and discuss that the 
compact supported perturbations used in \cite{lin} can be
achieved by Lagrangian perturbations from symplectic topology
point of view. Hence we get that the Casson-Lin's invariant is
Casson-Walker type invariant with non correction terms.

The only singular strata with trace zero condition consist of
{\re}s with image in the $U(1)$ subgroup and trace free.
Let ${\CS}^*(X)^{[i]} = Hom(\pi_1(X), U(1))^{[i]}/Z_2$ be the $U(1)$ strata of
$X$. Similar to (\ref{irre}), we have the following diagram, 
\begin{equation} \label{redu}
\begin{array}{ccc}
{\CS}^*(S^3 \setminus K)^{[i]} & \to & {\CS}^*(D_+^3 \setminus K)^{[i]} \\
\downarrow & & \downarrow \\
{\CS}^*(D_-^3 \setminus K)^{[i]} & \to & {\CS}^*(S^2 \setminus K)^{[i]} .
\end{array} 
\end{equation}
From Lin's notation $S_n = {\CS}^*(S^2 \setminus K)^{[i]}$ in \cite{lin},
the subset of $(X_1, \cdots, X_n, Y_1, \cdots, Y_n) \in H_n$ 
such that there is a matrix
$A \in SU(2)$ with property
\[A^{-1}X_iA, \ \ A^{-1}Y_iA \in U(1), \ \ 
\mbox{for all $i = 1, 2, \cdots, n$}. \]
Note that $Hom(\pi_1(S^3 \setminus K), U(1)) \cong U(1)$, 
${\CS}^*(S^3 \setminus K)^{[i]} = U(1)/Z_2$ the upper half unit circle in
the complex plane. There is 
only one reducible conjugacy class of {\re}s 
$\r : \pi_1(S^3 \setminus K) \to U(1) \hookrightarrow SU(2)$
which satisfies the
trace condition $tr (\r ([m_{x_i}])) = tr (\r ([m_{y_i}])) = 0$. 
Without loss of
generality, we may assume that $\r : \pi_1(S^2 \setminus K) \to U(1)$ is the
diagonal matrix
\[ \r ([m_{x_i}])  = \r ([m_{y_i}]) = \left[ \begin{array}{cc}
i & 0 \\
0 & - i  \end{array} \right]  .\]
Write $\r = \s \oplus \s^{-1}$ where $\s : \pi_1(S^2 \setminus K) \to U(1)$
is given by $\s ([m_{x_i}]) = \s ([m_{y_i}]) = i$. Then
\[ ad \r = End^0(\s \oplus \s^{-1}) = {\R } \oplus \s^{\ox 2} .\]
From deformation theory, any deformation in $H^1 (\cdot , {\R})$ 
direction changes the trace condition.
So only $H^1 (\cdot , \s^{\ox 2})$ fit into the Mayer-Vietoris sequence:
\[ H^1(S^3 \setminus K, \s^{\ox 2}) \to H^1(D^3_- \setminus K, \s^{\ox 2})
\oplus
H^1(D^3_+ \setminus K, \s^{\ox 2}) \to
H^1(S^2  \setminus K, \s^{\ox 2}) \to  \cdots . \]
That ${\CS}^*(D_+^3 \setminus K)^{[i]}$ and ${\CS}^*(D_-^3 \setminus K)^{[i]}$ 
intersect transversally in $T_{\r}{\CS}^*(S^2 \setminus K)^{[i]}$ is equivalent 
to $H^1(S^3 \setminus K, \s^{\ox 2}) = 0$. According to Milnor \cite{mil},
$H^1(S^3 \setminus K, \s^{\ox 2})$ can be computed by considering the
cohomology of the infinite cyclic covering space $\widehat{S^3 \setminus
K}$. This last cohomology $H^1(\widehat{S^3 \setminus K})$ is a module of the
Laurent polynomial ring ${\Z}[t, t^{-1}]$.
\[H^1(\widehat{S^3 \setminus K}) \ox_{{\Z}[t, t^{-1}]} C =
H^1(S^3 \setminus K, \s^{\ox 2}) . \]
Therefore $0 = H^1(\widehat{S^3 \setminus K}) = {\Z}[t, t^{-1}]/{\D }_K(t)$
if and only if $\D_K(-1) \neq 0$ which is true for any knot $K$.
\begin{lm} [Lin \cite{lin} Lemma 1.4] \label{single}
For any $K = \overline{\b }$, 
${\CS}^*(D_+^3 \setminus K)^{[i]}$ and ${\CS}^*(D_-^3 \setminus K)^{[i]}$
intersect transversally in $T_s{\CS}^*(S^2 \setminus K)^{[i]}$. 
${\CS}^*(S^3 \setminus K)^{[i]} = \{s\}$ is an isolated 
point of the intersection
between ${\CS}^*(D_+^3 \setminus K)^{[i]}$ and 
${\CS}^*(D_-^3 \setminus K)^{[i]}$ in ${\CS}^*(S^2 \setminus K)^{[i]}$. 
\end{lm}

\noindent{\bf Remark:} The isolated point 
$\{s\} = {\CS}^*(S^3 \setminus K)^{[i]}$ plays the same role as the trivial 
connection in the Casson's invariant for integral homology 3-spheres. But
this $U(1)$ reducible may have special isotopy around $s$ which reflects
global twisting in the normal bundle of ${\CS}^*(S^2 \setminus K)^{[i]}$
as in rational homology 3-spheres situation \cite{wa}.

The vector spaces $\{H^1(S^2 \setminus K, \s^{\ox 2}_p)| p \in 
{\CS}^*(S^2 \setminus K)^{[i]}\}$ form a symplectic vector bundle
$\nu $ over ${\CS}^*(S^2 \setminus K)^{[i]}$. There is a Hermitian structure 
on $\nu $ compatible with its symplectic structure $\o $.
We call an isotopy $\{h_t\}_{0\leq t \leq 1}$ of ${\CR}(S^2 \setminus K)^{[i]}$
to be a {\bf special} if $h_t|_{{\CS}(S^2 \setminus K)^{[i]}} = Id$ 
for all $t$, and
in a neighborhood of ${\CS}^*(S^2 \setminus K)^{[i]}$, 
$h_t$ is induced by a complex
symplectic bundle automorphism of $\nu$. 
Proposition 1.20 in \cite{wa} shows that there exists a special
isotopies $\{h_t\}_{0\leq t \leq 1}$ of ${\CR}(S^2 \setminus K)^{[i]}$
such that $h_1({\CR}^*(D^3_+ \setminus K)^{[i]})$ is transverse to
${\CR}^*(D^3_- \setminus K)^{[i]}$.

We are going to 
discuss the proper perturbation around the reducible $s$ in order to obtain
the computation of $\lam_{CL}(\b)$.
To obtain the surgery formula for $h(\b)$, one can consider the expression
\begin{equation} \label{88}
 \lc ( {\a }^{-1} \b) - \lc (\b) = \# \tilde{\G }_{{\a }^{-1} \b} \cap
\overline{\G}_{id} - \# \tilde{\G }_{\b} \cap \overline{\G}_{id} ,
\end{equation}
where ${\a }^{-1} = \s_1^2 $ a full twist on the first two strands. More
generally, we can consider $\lc ( {\a }^{-1} \b) - \lc (\b)$ without
$\a = \s_1^{-2}$. However, for $\overline{{\a }^{-1} \b}$ to represent a knot,
we must have $\a$ lying in the pure braid group $\ker ( \pi : B_n \to S_n)$. 
Since
$\ker (\pi : B_n \to S_n)$ is generated by 
$\s_1^2 $ and its conjugate, it is enough
to get a formula for $\lc (\s_1^2 \b) - \lc (\b )$ 
where $\s_1^2 \b$ is obtained from
$\b$ by doing surgery along a loop circulating the first two strands.
\medskip

Given $\a \in B_n$, we have an induced automorphism on $\overline{H}_n$ defined
by
${\a }_* : \overline{H}_n \to \overline{H}_n$:
\[{\a }_* (X_1 \cdots X_n, Y_1 \cdots Y_n) = (X_1 \cdots X_n, \a (Y_1) \cdots
\a (Y_n) ). \]
Under this automorphism, we have
${\a }_*( \overline{\G}_{id} \cap \overline{\G }_{{\a }^{-1} \b}) =
\overline{\G}_{\a } \cap \overline{\G}_{\b}.$
As we perturb $\overline{\G }_{{\a }^{-1} \b}$ to a transverse position
$\tilde{\G }_{{\a }^{-1} \b}$, its image ${\a }_*(\tilde{\G }_{{\a }^{-1} \b})$
becomes transverse to $\overline{\G}_{\a }$. In addition, near the neighborhood
of reducible point $s$, the counter clockwise motion in perturbing
$\tilde{\G }_{{\a }^{-1} \b}$ is also preserved by ${\a }_*$. Hence
\begin{equation} \label{move}
 \# (\overline{\G}_{id} \cap  \tilde{\G }_{{\a }^{-1} \b}) =
\# (\overline{\G}_{\a } \cap \tilde{\G}_{\b}) .
\end{equation} 
Using (\ref{move}), we can rewrite (\ref{88}) as
\begin{eqnarray*}
\lc ( {\a }^{-1} \b) - \lc (\b) & = & 
\# (\overline{\G}_{id} \cap  \tilde{\G }_{{\a }^
{-1} \b}) - \# (\overline{\G}_{id} \cap \tilde{\G}_{\b}) \\
& = & \# (\overline{\G}_{\a } - \overline{\G}_{id}) \cap \tilde{\G}_{\b}.
\end{eqnarray*}

In some sense, $\overline{\G}_{\a } - \overline{\G}_{id}$
can be thought of as a cycle except it contains the reducible point $\r $ and
so is $\tilde{\G}_{\b}$. To remedy this situation, we consider the difference:
\begin{eqnarray*}
[\lc ( {\a }^{-1} \b) - \lc (\b)] - 
[\lc ( {\a }^{-1} \b^{'}) - \lc (\b^{'})] & = &
\# (\overline{\G}_{\a } - \overline{\G}_{id}) \cap \tilde{\G}_{\b}
 - \# (\overline{\G}_{\a } - \overline{\G}_{id}) \cap \tilde{\G}_{\b^{'}} \\
& = & \# (\overline{\G}_{\a } - \overline{\G}_{id}) \cap
(\tilde{\G}_{\b} - \tilde{\G}_{\b^{'}}).
\end{eqnarray*}
 
At a neighborhood of the reducible point $s$, we choose a
1-parameter family of $\tilde{\G}_{\b}(t)$ which brings
$(\tilde{\G}_{\b^{'}})_{\r }$ to $(\tilde{\G}_{\b})_{\r }$. Using this isotopy,
we can construct a cycle $(\tilde{\G}_{\b} - \tilde{\G}_{\b^{'}})$ globally
by taking $\tilde{\G}_{\b} - \tilde{\G}_{\b^{'}}$ outside a neighborhood of
$s$ and in the neighborhood connecting up $\bd \tilde{\G}_{\b}$ to
$\bd \tilde{\G}_{\b^{'}}$, by the isotopy $\bd \tilde{\G}_{\b}(t)$.
 
As $(\tilde{\G}_{\b} - \tilde{\G}_{\b^{'}})_{global}$ is a cycle in the
nonsingular part $\overline{H}_n$, the intersection number
\[ I(s) = \# (\overline{\G}_{\a } - \overline{\G}_{id})
\cap (\tilde{\G}_{\b} - \tilde{\G}_{\b^{'}})_{global} , \]
makes sense as the intersection as a relative cohomology class with
an absolute class. $I(s)$ also counts irreducibles via
the global twisting around the reducible
$s$ which is essentially the Walker's correction term. Although the Walker's
geometric construction of his correction term is not defined 
for this situation (the trace zero {\re}s contains non-typical element such as
trivial {\re} for rational homology 3-spheres),
we can use the analytic definition defined in \cite{clm} for the
Walker's correction term $I(s)$.
\begin{equation} \label{clw}
[\lc ( {\a }^{-1} \b) - \lc (\b)] - 
[\lc ( {\a }^{-1} \b^{'}) - \lc (\b^{'})] 
= \# (\overline{\G}_{\a } - \overline{\G}_{id})
\cap (\tilde{\G}_{\b} - \tilde{\G}_{\b^{'}})_0 + I(s).
\end{equation}

\begin{enumerate}
\item $\# (\overline{\G}_{\a } - \overline{\G}_{id})
\cap (\tilde{\G}_{\b} - \tilde{\G}_{\b^{'}})_0$ can be expressed in turn of
triple Maslov index
\begin{equation} \label{89}
\# (\overline{\G}_{id} - \overline{\G }_{\a}) \cap (\tilde{\G}_{\b} -
\tilde{\G}_{\b^{'}})_0 = 
\frac{1}{2} Mas(\overline{\G }_{\a}, \tilde{\G}_{\b}, \overline{\G}_{id})_{s}
-  \frac{1}{2} Mas(\overline{\G }_{\a},
\tilde{\G}_{\b^{'}}, \overline{\G}_{id})_{s} .
\end{equation}
This is because, from the definition, $\tilde{\G}_{\b}, \tilde{\G}_{\b^{'}}$,
are transverse to $\overline{\G }_{\a}, \overline{\G}_{id}$.
Hence the dimension correction terms all disappear 
(see \cite{clm} and \cite{clm1}).
\item $\lc (\a^{-1} \b^{'}) - \lc (\b^{'})$
in (\ref{clw}) can be eliminated by choosing ${\b }^{'}$ as follows.
There exists a braid $\b^{'} $ such that both
$\overline{\b^{'}}$ and $\overline{\s^2_1 \b^{'}}$ represent the trivial knot.
In the case $n =2$,
we have choose $\overline{\b^{'}}$ and $\overline{\s^2_1 \b^{'}}$ to be
$\overline{\s^{-1}_1}, \s_1$ respectively. For $n > 2$, we have 
$\b^{'} = \s^{-1}_1 \s^{-1}_2 \cdots \s^{-1}_n$,
$\overline{\s^2_1 \b^{'}} = \s_1 \s^{-1}_2 \cdots \s^{-1}_n$.
Thus $\lc (\a^{-1} \b^{'}) - \lc (\b^{'}) = 0$ and $Mas(\overline{\G }_{\a},
\tilde{\G}_{\b^{'}}, \overline{\G}_{id})_{s} = 0$.
(\ref{clw}) is reduced to 
\[  \lc (\a^{-1} \b) - \lc (\b) = 
\frac{1}{2} 
Mas(\overline{\G }_{\a}, \tilde{\G}_{\b}, \overline{\G}_{id})_{s} + I(s) .\]
\item Using Wall's nonadditivity of signature, 
Cappell, Lee and Miller in 
\cite{clm} identified the triple Maslov index with the (twisted) signature.
So (\ref{clw}) becomes the following.
\begin{equation} \label{817}
\lc (\a^{-1} \b) - \lc (\b)  = 
\frac{1}{2} \{ sign (\overline{\s^2_1 \b}) -
sign(\overline{\b}) \} + I(s).
\end{equation}
\item $I(s) = 0$ and the signature of $\overline{\s_1^2}$
is clear zero.
\[ \lc (\b) = \frac{1}{2} sign (\overline{\b}; [i^2]) 
= \frac{1}{2} sign (\overline{\b}) . \]
This is Corollary 2.10 in \cite{lin}. 
In the following, we show that the Walker-type 
correction term $I(s)$ vanishes in two different methods.
\end{enumerate}

\begin{lm} \label{812}
For $\a = \s^{-2}_1$, $I(s) = 0$. (see also \cite{clm})
\end{lm}
Proof: For $\a = \s^{-2}_1$, we have
\begin{eqnarray*}
\a (X_1) & = & (X_1X_2)^{-1} X_1 (X_1X_2) \\
\a (X_2) & = & (X_1X_2)^{-1} X_2 (X_1X_2) \\
\a (X_j) & = & X_j, \ \ \ j \geq 3 .
\end{eqnarray*}
There is a Casson-deformation $X^t$ from $SU(2) \setminus \{-I \}$ to itself
such that $X^t$ commutes with $X$, 
\begin{eqnarray*}
\p (t)(X_1) & = & [(X_1X_2)^t]^{-1} X_1 [(X_1X_2)^t] \\
\p (t)(X_2) & = & [(X_1X_2)^t]^{-1}X_2 [(X_1X_2)^t] \\
\p (t)(X_j) & = & X_j, \ \ \ j \geq 3 .
\end{eqnarray*}
Since trace $tr(\p (t)(X_j)) = tr(X_j)$ for all $j$,
the formula
\[\overline{\G}_{\p (t)} = \{(X_1, \cdots, X_n,\p (t)(X_1), \cdots,
\p (t)(X_n)) \in Q_n \x Q_n\}/SU(2) , \]
gives us a 1-parameter family of subspaces which connects up ${\G}_{id}/SU(2)$
at $t = 0$ to ${\G}_{\a}/SU(2)$ at $t = 1$. As 
$\p (t)(X_1X_2) = [(X_1X_2)^t]^{-1}(X_1X_2) [(X_1X_2)^t] = X_1X_2,
\p (t)(X_j) = X_j, j \geq 3$
we have
\[\p (t)(X_1X_2\cdots X_n) = (X_1X_2\cdots X_n) , \]
and so ${\G}_{\p (t)}$ is a subspace in ${\CR}(S^2 \setminus K)^{[i]}$.
By definition of the Casson-deformation, 
$$\p (t)|_{{\CS}^*(S^2 \setminus K)^{[i]}} = id, $$
so $\p (t)$ is a special isotopy.
It follows that
$\{{\G }_{\p (t)} : 0 \leq t \leq 1\}$
gives us a cycle in ${\CR}(S^2 \setminus K)^{[i]}$ whose boundary is 
${\G}_{id}/SU(2) - {\G }_{\a}/SU(2) = {\CR}(D^3_- \setminus K)^{[i]}
- {\G }_{\a}/SU(2)$. Thus
\begin{eqnarray*}
- I(s) & = & \# ({\CR}(D^3_- \setminus K)^{[i]} - {\G }_{\a}/SU(2))
\cap (\tilde{\G}_{\b} - \tilde{\G}_{\b^{'}})_{global} \\
 & = & \# (\bd \{{\G}_{\p (t)}: 0 \leq t \leq 1\} \cap
(\tilde{\G}_{\b} - \tilde{\G}_{\b^{'}})_{global} \\
 & = & 0 .
\end{eqnarray*}
\qed

Using spectral flows in \cite{clm1}, the Walker correction term $I(s)$ for
$s = \s \oplus \s^{-1}$, can be computed as the rho invariant
$\r (M, \s^{\otimes 2})$ (see Proposition 7.1 in \cite{clm1} Part III). 
\begin{thm} [Cappell, Lee and Miller \cite{clm1}] \label{rho}
Let $s = \s \oplus \s^{-1}$ be a $U(1)$ {\re} of $\pi_1(M)$ for the
rational homology 3-sphere $M$. Then the Walker's correction term
$I(s)$ is given by
\[ I(s) = - \r (M, \s^{\otimes 2}) , \]
where $\r (M, \s^{\otimes 2})$ is the rho invariant of $M$ associated
to the {\re} $\s^{\otimes 2}$.
\end{thm}

We take the above as analytic definition for the Walker-type
invariant. Applying Theorem~\ref{rho} to $M = \hat{S^3}$ the double
branched cover of $S^3$ with branched set $K$ and $s = \s \oplus \s^{-1}$ 
which is the fixed point of $Z_2$ operation on ${\CS}(\hat{S^3})^{Z_2}$, 
we have that $ad s = {\R} \oplus \s^{\otimes 2}$ acts by 1 on the
real part of ${\bf su(2)}$ and acts by
$\s^{\otimes 2} = -1$ on the complex part of ${\bf su(2)}$.
Thus for $M = \bd X^4$, $\s^{\otimes 2} = -1$ is trivially extended to a
unitary {\re} of $\pi_1(X)$. 
By Theorem 2.4 \cite{aps}, 
\begin{eqnarray*}
I (s) & = & - \r (M, \s^2) \\
& = & - (sign_{Id} (X) - sign_{\s^{\ox 2}} (X)) \\
 & = & 0 .
\end{eqnarray*} 
Note that the signature of Hermitian form $x^* {\G}_z y$ with ${\G}_z =
\frac{1}{2} \{ (1 - \overline{z}) \G + (1 -z) \G\}$ has the relation
$sign_{\xi} = sign(1 - Re \xi) sign_z$ for 
$z = \frac{1 -\xi}{1- \overline{\xi}}$, the only interesting case is when 
$|z| = 1$, where $\G_z$ reduces to $\frac{1}{2} (1 - \overline{z}) \G$. In 
particular for $z = -1$, the signature is the same as $sign_{Id}$. See
\cite{mat} for example.

Let $\r \in {\CS}^*(S^2 \setminus K)^{[i]}$ be the $U(1)$ reducible
{\re} of $\pi_1(S^2 \setminus K)$ with trace zero condition.
A normal neighborhood of $\rho $ can be isomorphic to 
the cone of $H^1(S^2 \setminus K, {\bf h}^{\perp}_{Ad \rho})/U(1)$,
(${\bf h}$ is the Lie algebra of the fixed oriented maximal torus of $SU(2)$,
${\bf h}^{\perp}$ is the orthogonal complement of ${\bf h}$ with respect to
the Killing form of $SU(2)$),  
the cone bundle $E({\CS}^*(S^2 \setminus K)^{[i]}) \to 
{\CS}^*(S^2 \setminus K)^{[i]}$ is isomorphic to a neighborhood of 
${\CS}^*(S^2 \setminus K)^{[i]}$ in ${\CR}(S^2 \setminus K)^{[i]}$
via an exponential map $\exp : {\cal N}({\CS}^*(S^2 \setminus K)^{[i]})
\to Q_n \x Q_n$. 
The exponential map allows us to identify a germ of functions $\phi
\in C^{\infty }(Q_n \x Q_n, {\R})$ near 
${\CS}^*(S^2 \setminus K)^{[i]}$ with a function on 
${\cal N}({\CS}^*(S^2 \setminus K)^{[i]})$. For example,
choose a partition of unity $\c : Q_n \x Q_n \to {\R}$
which is $0$ outside ${\cal N}({\CS}^*(S^2 \setminus K)^{[i]})$ and 
$1$ near the 0-section
and consider a function $g: {\cal N}({\CS}^*(S^2 \setminus K)^{[i]}) \to {\R}$
which is induced by a
Hermitian pairing on each fiber 
$H^1(S^2 \setminus K, {\bf h}^{\perp}_{Ad \rho})/U(1)$. 
Locally, $g(\rho, z) = \sum a_{ij}(\rho) z_i \overline{z_j}, \ \ z=(z_1, \cdots,
z_d) \in H^1(S^2 \setminus K, {\bf h}^{\perp}_{Ad \rho})/U(1)$,
\begin{equation} \label{qua}
H_g (\exp (\rho, z)) = \c (\rho ,z) g(\rho , z)
\end{equation}
gives us a function $H_g \in C^{\infty }(Q_n \x Q_n, {\R})$.
The vector field $grad H_g$ associated to
this perturbation has the property that $g(\rho , z)$ 
is quadratic in the normal $z$-direction and so $grad H_g
= 0$ when it is restricted to the zero section, the
Hessian Hess$(H_g)$ of $H_g$ at each fiber 
$H^1(S^2 \setminus K, {\bf h}^{\perp}_{Ad \rho})/U(1)$
is given by the Hermitian pairing $\sum a_{ij}(\rho) z_i
\overline{z_j}$.

Let ${\cal L}ag(W)$ be the space of complex Lagrangians
in $W = H^1(S^2 \setminus K, {\bf h}^{\perp}_{Ad \rho})$.
${\cal L}ag(W)$ is a homogeneous manifold
whose tangent space can be identified with the space of Hermitian pairings
on $H^1(S^2 \setminus K, {\bf h}^{\perp}_{Ad \rho})$
(c.f. \cite{wa} p13-14).
We consider the open subspace $\Pi_0$ of perturbation
data
in $\Pi $ such that
\begin{enumerate}
\item $grad H_g( \r ) = 0$ for $\r  \in {\CS}^*(S^2 \setminus K)^{[i]}$,
\item $(\p_{\b} + H_g)(\overline{\G}_{id})$ is mapped 
injectively into a complex
Lagrangian subspace in $H^1(S^2 \setminus K, {\bf h}^{\perp}_{Ad
\rho})$.
\end{enumerate}
Let $\mbox{Lag}({\cal W})$ denote the smooth
fiber bundle over ${\CS}^*(S^2 \setminus K)^{[i]}$ whose fiber at
a point $\r \in {\CS}^*(S^2 \setminus K)^{[i]}$ 
is the homogeneous space ${\cal L}ag(W)$, and
$\Gamma(\mbox{Lag}({\cal W}))$
denote the space of smooth sections of this bundle completed with respect
to an appropriate Sobolev norm. Then by assigning to $\pi$ the section of
Lagrangian $\r \mapsto (\p_{\b} + H_g)(\overline{\G}_{id})$ in 
$H^1(S^2 \setminus K, {\bf
h}^{\perp}_{Ad
\rho})$ there exists a mapping
$p : \Pi_0 \to \Gamma(\mbox{Lag}({\cal W}))$.
\begin{pro} \label{ind}
(i) The mapping $p : \Pi_0 \to \Gamma(\mbox{Lag}({\cal W}))$
is a submersion whose image is the entire
space $\Gamma(\mbox{Lag}({\cal W}))$.

(ii) $\Pi_0$ is path connected, i.e. any two different perturbations can be 
homotopically connected to each other through a 1-parameter family in 
$\Pi_0$.
\end{pro}
Proof: The argument  follows exactly as in the proof of Theorem B and 
Proposition 3.3 in \cite{ll}. Since the correction term $I(s)$ 
vanishes, there is non obstruction to connect any two different perturbations.
\qed

Now the Walker-type correction term vanishing implies that 
the compact support perturbation $H_g$ does not effect the invariant 
$\lc (\b)$. Also the Floer homology defined in the next section will
be independent of the Lagrangian perturbation around the reducible
$s$, unlike the situation studied in \cite{ll}.
So we can choose the counter clock orientation (pick $H_g$ the
quadratic form in the positive direction) to perturb 
${\CR}^*(D^3_+ \setminus K)^{[i]}$ (as graph of $\p_{\b}$)
into $\hat{{\CR}}^*(D^3_+ \setminus K)^{[i]}$
(the graph of $\ti{\p}_{\b}$) which meets 
${\CR}^*(D^3_- \setminus K)^{[i]}$ transversely. Hence
the Casson-Lin's invariant $\lc (\b)$ is also the Casson-Walker type invariant
with trivial correction term $I(s)$.

\section{The symplectic Floer homology of a braid}
 
\subsection{Floer homology for symplectic fixed points}

Let $H: {\R} \x M \to \R$, $H(s, x) = H_s(x)$ be a smooth time-dependent
Hamiltonian function such that $H_s = H_{s+1} \circ \p$. The symplectic
diffeomorphism $\psi_s : M \to M $ generated by $H$ are defined by
\[\frac{d \psi_s}{ds} = X_{H,s} \circ \psi_s, \]
\[ \psi_0 = Id, \ \ \ \o (X_{H,s}, \cdot ) = dH_s . \]
They satisfy $\psi_{s+1} \circ \p_H = \p \circ \psi_s$, where
$\p_H = \psi_1^{-1} \circ \p$. For a generic Hamiltonian $H_s$, the fixed points
of $\p_H$ are all nondegenerate (\cite{fl3}, \cite{hs}), denoted by $Fix(\p_H)$.
They can be identified with the critical points of the perturbed symplectic
action functional on the space of smooth paths
\[\O_{\p} = \{ \g : {\R} \to M \ \ | \ \ 
\g (s+1) = \p (\g (s)) \} .\]
If $M$ is simply connected, then $\pi_1 (\O_{\p} ) = \pi_2 (M)$; otherwise, we
will apply the construction to components of contractible loops of $M$.
The perturbed symplectic action functional $a_H: \O_{\p} \to {\R}/\a {\Z}$ is
defined as a function whose differential is 
\[ da_H(\g) \xi = \int_0^1 \o (\dot{\g} - X_s(\g), \xi) ds. \]
So the critical points of $a_H$ are the paths of the form $x(s) = \psi_s (x_0)$
such that $x(s+1) = \p (x(s))$ which are in one-to-one correspondence with
the fixed point set $Fix(\p_H)$. 

Now choose a smooth map ${\R} \to {\J}(M, \o): s \mapsto J_s$ such that
$J_s = \p^* J_{s+1}$. Such a structure determines a $L^2$-metric on $\O_{\p}$
and the gradient of $a_H$. In particular, the gradient flow of $a_H$ can be 
identified with the solution $u: {\R}^2 \to M$ of the pseudoholomorphic curves
\begin{equation} \label{pde}
\pjh (u) = \bd_t u + J_s(u) (\bd_s u - X_s (u)) = 0, 
\end{equation}
\begin{equation} \label{bc}
u(s+1, t) = \p (u (s,t)) , 
\end{equation}
in the sense of Gromov (see \cite{gr} for $X_s = 0$ case and \cite{fl3} for
general case). If all fixed points in $Fix(\p_H)$ are nondegenerate, then any 
solution of (\ref{pde}) and (\ref{bc}) with finite energy
\begin{equation} \label{fe}
E(u) = \frac{1}{2} \int_{{\R}} \int^1_0 (|(\bd_s u - X_s (u))|^2 +
|\bd_t u|^2) ds dt  < \infty 
\end{equation}
has limit, by Gromov's compactness (see \cite{fl3} and \cite{gr}), 
\begin{equation} \label{lv}
\lim_{t \to \pm \infty} u(s,t) = \psi_s (x^{\pm }), \ \ \
x^{\pm } = \p_H (x^{\pm }),  
\end{equation}

Let ${\CM}(x,y)$ be the 
space of solutions $u$ of (\ref{pde}), (\ref{bc}) and (\ref{fe}) with
$x, y \in Fix(\p_H)$. This gives a perturbed Cauchy-Riemann operator
$D_u: L^p_{1, \p}(u^*TM) \to L^p_{0, \p}(u^*TM)$ defined by
\begin{equation} \label{cro}
D_u \xi = \n_t \xi + J_s(u) (\n_s \xi - \n_{\xi }X_s (u)) + 
\n_{\xi } J_s (u) (\bd_s u - X_s (u)), 
\end{equation}
where $\n$ denotes the covariant derivative with respect to the $s$-dependent
metric $<v_1, v_2>_s = \o (v_1, J_s v_2)$, and $L^p_{0, \p}(u^*TM)$
(respectively $L^p_{1, \p}(u^*TM)$) is the completions of the space of smooth
vector fields $\xi (s,t) \in T_{u(s,t)}M$ along $u$ with $\xi (s+1, t) =
d\p (u(s,t)) \xi(s,t)$ and compact support on $S^1 \x {\R}$, with respect to
the $L^p$-norm (respectively $L^p_1$-norm) on $S^1 \x {\R}$.
If $x, y \in Fix(\p_H)$ are nondegenerate and $u$ satisfy (\ref{pde}), (\ref{bc})
and (\ref{fe}), then $D_u$ is Fredholm and its index is given by the Maslov
class of $u$:
\[ index D_u = \mu_u (x, y) .\]
The Maslov class $\mu_u$ is invariant under homotopy, additive for catenations 
and satisfies
\begin{equation} \label{mcn}
\mu_{u \# v} = \mu_u - 2 c_1(v) , 
\end{equation} 
for $v \in \pi_2 (M)$ by \cite{fl3} Proposition 1b. A Hamiltonian function
$H$ is called regular if the fixed point of $\p_H$ are all nondegenerate
and the operator $D_u$ is onto for every $u \in {\CM}(x, y)$ and $x, y \in
Fix(\p_H)$. Floer in \cite{fl3} showed that the
set ${\cal H}^{reg} = {\cal H}^{reg}(J)$ of regular Hamiltonian functions
is generic in the sense of Baire with respect to a 
suitable $C^{\infty}_{\ve}$-topology. 
For $H \in {\cal H}^{reg}$, $x, y \in Fix(\p_H)$, the space ${\CM}(x, y)$ is
a manifold whose local dimension near $u$ is the Maslov class
$\mu_u(x, y)$. By (\ref{mcn}), 
\[\mu : Fix(\p_H) \to {\Z}_{2N}, \]
defined up to an additive constant, such that
\[ \mu_u (x, y) = \mu (x) - \mu (y) \ \ \ \pmod {2N} , \]
for every $u \in {\CM}(x, y)$. Here the integer $N$ is the minimal Chern
number defined by
$I_{c_1}: \pi_2(M) \to {\Z}$, $Im I_{c_1} = <N>$ since $\Z$  is PID and every
ideal is generated by an element in $\Z$.

Let $C_i$ be the free generated $\Z$ module over $x \in Fix(\p_H)$ and
$\mu (x) \equiv i \pmod {2N}$. Then the Floer boundary map
\begin{equation} \label{fbc}
\bd : C_i \to C_{i-1}, \ \ \ x \mapsto \sum_{\mu(x) - \mu (y) =1}
\# \hat{{\CM}}(x,y) y ,
\end{equation}
is defined by taking the sum of the algebraic number over all
one dimensional components of $\hat{{\CM}}(x, y)$, 
where $\hat{{\CM}}(x, y)$ is the balanced moduli space of $J$-holomorphic
curves, see also \cite{fl3} and \cite{ll}. These numbers are defined
by comparing the
flow orientation of $u$ with the coherent orientation of ${\CM}(x, y)$
as in \cite{fh}. One can check that $\bd $ is well-defined
and $\bd \circ \bd =0$ as Floer did in \cite{fl3} for $\p = id$ case.
The homology of this chain complex $(C_i, \bd )$ is defined to be the Floer
homology for symplectic fixed points. I.e.
\[HF^{sym}_* (M, \p, H, J) = \ker \bd / Im \bd , \ \ \ * \in Z_{2N} .\]
The Floer homology groups are independent of the almost complex structures
$J_s$ and the perturbation $H \in {\cal H}^{reg}$ used in the construction
\cite{fl3}. They depend on $\p$ only up to Hamiltonian isotopy. There
is a natural isomorphism
\[ HF^{sym}_* (M, \p^0, H^0, J^0) = HF^{sym}_* (M, \p^1, H^1, J^1) ,\]
whenever $\p^0$ and $\p^1$ are related by a Hamiltonian isotopy. 

For the monotone symplectic manifold $(M, \o)$, there are $2N$-graded
Floer homology groups $HF_*^{\mbox{sym}}(\p )$ for every symplectic
diffeomorphism $\p$, where $N$ is the minimal Chern number. The Euler
characteristic of the Floer homology is the Lefschetz number of $\p$,
$$\chi (HF^{sym}_* (M, \p)) = L(\p ). $$
If there are degenerated fixed points, then one can choose a Hamiltonian 
perturbation.  
For any two symplectic diffeomorphisms $\p , \psi$, 
there is a natural isomorphism
\begin{equation} \label{qdc}
HF_*^{\mbox{sym}}(\p ) \stackrel{\cong}{\rightarrow} 
HF_*^{\mbox{sym}}(\psi \circ \p \circ \psi^{-1} ) .
\end{equation}
(see \cite{ms} for the so called Donaldson quantum category)

Note that the symplectic Floer homology is essentially an infinite dimensional
version of Morse-Novikov theory for the symplectic action functional. The unique
covering space $\ti{\O}_{\p}$ of the space $\O_{\p}$ of the contractible loops 
has deck-transformation group $\G = Im (\pi_2(M) \to H_2(M))$. The symplectic
action functional is well-defined over $\ti{\O}_{\p}$, then applying the
Novikov's construction, one obtains the symplectic Floer homology as a module
over the Novikov ring $\Lam_{\o}(\G)$ (see \cite{ms} and \cite{vo}).
For our case, we do not need to extend the symplectic Floer homology
with Novikov ring coefficients.                      

\subsection{Floer homology of braids}

In this subsection, we are going to apply the construction in \S 4.1 to 
the monotone symplectic manifold ${\CR}^*(S^2 \setminus K)^{[i]}$ and
the symplectic diffeomorphism $\p_{\b}$ induced from a braid $\b$.

For any knot $K = \overline{\b}$, $\b \in B_n$ a braid, we apply the previous
section to $\O_{\p_{\b}}$ the space of contractible loops in
${\CR}^*(S^2 \setminus K)^{[i]}$. If there are degenerated fixed points of
$\p_{\b}$, then we apply the Lagrangian perturbation around $s$, as 
discussed in \S 3, it gives arise a Hamiltonian vector field perturbation 
to make all the fixed points of $\p_{\b}$ to be nondegenerated.
Instead of relative index with shift, we can assign an relative index 
for each element in Fix$(\ti{\p}_{\b})$ with respect to the
$U(1)$-{\re} $s$, $\mu (x) = \mu_u (x, s) \pmod {2N}$. Thus the symplectic
Floer chain complex is
\[C_i = \{ x \in Fix(\ti{\p}_{\b}) \cap {\CR}^*(S^2 \setminus K)^{[i]} :
\mu (x) = i \in Z_{2N} \} . \]
\begin{pro} \label{wel}
There is a well-defined symplectic Floer homology $HF_*^{\mbox{sym}}(\p_{\b},
H, J), * \in Z_{2N}$ of a braid representing a knot.
\end{pro}
Proof: This basically follows from argument in  \cite{fl3} or \cite{ll1}.
We first show that $\bd$ is well-defined and does not involved the 
$U(1)$ reducible {\re} $s \in Fix(\ti{\p}_{\b})$. Recall the special isotopy
does fix the $U(1)$-strata and $\{s \}= {\CS}^*(S^3 \setminus K)^{[i]}$.
Because any trajectories connecting with $s$ always come in at least
two dimensional families. Followed the method in \cite{fl3} p 587, we can ensure
that these two dimensional families belong to regular Morse cells and their
index is at least $2$, so that $s$ does not contribute to $\bd $.
In order to show $\bd \circ \bd = 0$, we consider the term for a fixed
$y$ in the following.
\[ \bd \circ \bd (x) = \sum_{z \in C_{i-1}} \sum_{y \in C_{i-2}}
\# \hat{{\CM}}(x,z) \cdot \# \hat{{\CM}}(z, y) y .\]
For the pair $(x, y)$, there is the two dimensional moduli space of 
$J$-holomorphic curves ${\CM }^2(x, y)$. The ends of $\hat{{\CM}}(x,y)$ consists
of all the components $\hat{{\CM}}(x,z) \x \hat{{\CM}}(z,y)$
for $z \in Fix(\ti{\p}_{\b}) \cap {\CR}^*(S^2 \setminus K)^{[i]}$.
It is impossible for $z$ to be the reducible $s$ because the $U(1)$ symmetry
group would add to the extra one more parameter families 
inside the moduli space. Hence the standard argument shows that
$\bd \circ \bd = 0$. Another way to make sure that the reducible {\re} $s$ 
does not contribute is to apply the gluing technique developed in \cite{ll}
for the $U(1)$ reducible {\re}s. See \cite{ll} \S 4. \qed

In order to get an invariant of knots from braids, we have to verify that
$HF_*^{\mbox{sym}}(\p_{\b}, H, J)$ (independent of $(H,J)$) is invariant
under Markov moves. A Markov move of type I changes $\s \in B_n$ to 
$\xi^{-1} \s \xi \in B_n$ for any $\xi \in B_n$, and the Markov move of
type II changes $\s \in B_n$ to $\s_n^{\pm } \s \in B_{n+1}$, or the inverses 
of these operations. It is well-known that two braids $\b_1$ and $\b_2$ has 
isotopic closure if and only if $\b_1$ can be changed to $\b_2$ by a sequence 
of finitely many Markov moves \cite{bi}.
\begin{thm} \label{invar}
For $\overline{\b_1} = \overline{\b_2} = K$ as a knot, $\b_1 \in B_n, \b_2 \in
B_m$, then there is a natural isomorphism
\[ HF_i^{\mbox{sym}}(\p_{\b_1}) \cong HF_i^{\mbox{sym}}(\p_{\b_2}), \ \ \ 
i \in Z_{2N} .\]
So the symplectic Floer homology 
$\{HF_i^{\mbox{sym}}(\p_{\b})\}_{i \in Z_{2N}}$ is a knot invariant.
\end{thm}
Proof: We only need to show that for $\b \in B_n$ with $\overline{\b}$ being
a knot $K$, the Markov moves of type I and type II on $\b$ provide a Hamiltonian
isotopy of $\p_{\b}$. Hence from the invariance property of the symplectic 
Floer homology, we get that $\{HF_i^{\mbox{sym}}(\p_{\b})\}_{i \in Z_{2N}}$
is an invariant of knot $K = \overline{\b}$.

Suppose we have the Markov move of type I: change $\b$ to $\xi^{-1} \b \xi$
for some $\xi \in B_n$. The element $\xi$ in $B_n$ induces a diffeomorphism
$\xi: Q_n \to Q_n$ is orientation preserving as observed by Lin in \cite{lin}.
Note that $B_n$ is generated by $\s_1, \cdots, \s_{n-1}$. For any $\s_i^{\pm }$,
the induced diffeomorphism $\s_i^{\pm } \x \s_i^{\pm }: 
Q_n \x Q_n \to Q_n \x Q_n$ is an orientation preserving and 
symplectic diffeomorphism. So $\xi $ is also a symplectic orientation preserving
diffeomorphism since orientation preserving and symplectic form preserving
properties are invariant under the composition operation. Hence there is a
symplectic diffeomorphism
\[\xi \x \xi : Q_n \x Q_n \to Q_n \x Q_n , \]
which commutes with the $SU(2)$-action and
\[ \xi \x \xi ({\CR}^*(S^2 \setminus K)^{[i]}) = {\CR}^*(S^2 \setminus K)^{[i]}
 \ \ \ (\mbox{changing variables by $\xi \x \xi$}),\]
\[\xi \x \xi ({\CR}^*(D_-^3 \setminus K)^{[i]}) = 
{\CR}^*(D_-^3 \setminus K)^{[i]} \ \ \ 
(\mbox{in new coordinate $\xi(X_1), \cdots, \xi (X_n)$}),\]
\[ \xi \x \xi ({\CR}^*(D_+^3 \setminus K)^{[i]}) = 
{\CR}^*(D_+^3 \setminus K)^{[i]}\ \ \ 
(\mbox{in new coordinate $\xi(X_1), \cdots, \xi (X_n)$}),\]
as oriented manifolds. Let $f_{\xi} : {\CR}^*(S^2 \setminus K)^{[i]} \to 
{\CR}^*(S^2 \setminus K)^{[i]}$ be the induced symplectic diffeomorphism of 
$\xi \x \xi$. Hence we have
\[ \p_{\b} = f_{\xi} \circ \p_{\xi^{-1} \b \xi} \circ f_{\xi}^{-1} , \]
from changing variables via $f_{\xi}$. Note that $Fix(\p_{\xi^{-1} \b \xi})$
is identified with $Fix(\p_{\b})$ under $f_{\xi}$. So we get
$\ti{\p}_{\b} = f_{\xi} \circ \ti{\p}_{\xi^{-1} \b \xi} \circ f_{\xi}^{-1}$
if necessary for Hamiltonian perturbations. Therefore by (\ref{qdc}),
\begin{equation} \label{type1}
HF_i^{\mbox{sym}}(\p_{\xi^{-1} \b \xi}) \cong 
HF_i^{\mbox{sym}}(f_{\xi} \circ \p_{\xi^{-1} \b \xi} \circ f_{\xi}^{-1}) 
= HF_i^{\mbox{sym}}(\p_{\b}) .
\end{equation}
It is clear that the argument goes through for the inverse operation of Markov
move of type I.

Suppose we have the Markov move of type II: change $\b$ to 
$\s_n \b \in B_{n+1}$. Recall that $\s_n (x_i) = x_i, 1 \leq i \leq n-1,
\s_n (x_n) = x_n x_{n+1} x_n^{-1}$ and $\s_n (x_{n+1}) = x_n$. We need to 
identify the Floer homology from the construction in
$\hat{H}_n$ into the one from $\hat{H}_{n+1}$. Following Lin \cite{lin}, there
is an imbedding $g: Q_N \x Q_n \to Q_{n+1} \x Q_{n+1}$ given by
\[g(X_1, \cdots, X_n, Y_1, \cdots, Y_n) = (X_1, \cdots, X_n,Y_n,
Y_1, \cdots, Y_n, Y_n) .\]
Such an imbedding commutes with the $SU(2)$-action and $g(H_n) \subset H_{n+1}$,
and induces an imbedding
\[ \hat{g}: \hat{H}_n (={\CR}^*(S^2 \setminus \overline{\b})^{[i]}) \to
\hat{H}_{n+1} (= {\CR}^*(S^2 \setminus \overline{\s_n \b})^{[i]}) .\]
Note that the symplectic structure of $\hat{H}_{n+1}$ restricted on
$\hat{g}(\hat{H}_n)$ is the symplectic structure on $\hat{H}_n$. Hence
$\hat{g}$ is a symplectic imbedding which makes the space
$\hat{H}_n$ naturally into the symplectic submanifold $\hat{g}(\hat{H}_n)$
of $\hat{H}_{n+1}$. Under this imbedding, we have
$\hat{g}(\p_{\b}): \hat{H}_{n+1} \to \hat{H}_{n+1}$ is given by
\begin{equation} \label{pb}
(X_1, \cdots, X_n, X_1, \cdots, X_n) \mapsto (X_1, \cdots, X_n, \b(X_n),
\b(X_1), \cdots, \b(X_n), \b(X_n)).
\end{equation}
The image of $\hat{g}(\p_{\b})$ is invariant under the operation of $\s_n$.
Also the corresponding symplectic diffeomorphism $\p_{\s_n \b}$ is given by
\[\p_{\s_n \b}(X_1, \cdots, X_n, X_{n+1}, X_1, \cdots, X_n, X_{n+1}) \]
\begin{equation} \label{pnb}
 = (X_1, \cdots, X_{n+1}, \b(X_1), \cdots, \b(X_{n-1}), \b(X_n) X_{n+1} 
\b(X_n)^{-1}, \b(X_n) ) .
\end{equation}
Thus we have
\[ \hat{g}({\CR}^*(D_-^3  \setminus \overline{\b})^{[i]}) \subset
{\CR}^*(D_-^3  \setminus \overline{\s_n \b})^{[i]}, \ \ 
\hat{g}({\CR}^*(D_+^3  \setminus \overline{\b})^{[i]}) \subset
{\CR}^*(D_+^3  \setminus \overline{\s_n \b})^{[i]} . \]
The fixed points of $\p_{\s_n \b}$ are elements
\[ \b(X_i) = X_i, 1 \leq i \leq n_1; \ \ \ 
 \b(X_n) X_{n+1} \b(X_n)^{-1} = X_n, \ \ \ 
\b(X_n) = X_{n+1} , \]
which is equivalent to $\b(X_i) = X_i, 1\leq i \leq n$, i.e.
\[ Fix(\p_{\s_n \b}) = Fix (\hat{g}(\p_{\b})) = Fix (\p_{\b}) . \]
For degenerate fixed point of $\p_{\b}$, we perturb 
${\CR}^*(D_+^3  \setminus \overline{\b})^{[i]}$ via Hamiltonian vector field
with compact support in $\hat{H}_n$ so that all elements in
$Fix(\ti{\p}_{\b})$ are nondegenerated. By the standard isotopy 
extension argument, we can further perturb 
${\CR}^*(D_+^3  \setminus \overline{\s_n \b})^{[i]}$ ($\supset
\hat{g}({\CR}^*(D_+^3  \setminus \overline{\b})^{[i]})$) with compact
support in $\hat{H}_{n+1}$ such that
\[ \ti{{\CR}}^*(D_+^3  \setminus \overline{\s_n \b})^{[i]} \cap
\hat{g}(\hat{H}_n) = \hat{g} 
(\ti{{\CR}}^*(D_+^3  \setminus \overline{\b})^{[i]},\]
and all elements in $Fix(\ti{\p}_{\s_n \b})$ are nondegenerated. So we have 
\[Fix(\ti{\p}_{\s_n \b}) = Fix (\ti{\p}_{\b}) .\]
In \cite{lin}, Lin verified that orientations involved in the process are
preserved. Hence the algebraic number of $Fix(\ti{\p}_{\s_n \b})$ equals to 
one of $Fix(\ti{\p}_{\b})$. The corresponding relative Maslov indexes are
also unchanged under $\hat{g}$.

Take the standard 2-sphere $S^2 \subset {\R}^3$ and the height function
$H = x_3$ (determined a vector field $X_H$ by $\o (X_H, \cdot ) = d H$).
The level sets are circles at constant height and the Hamiltonian flow
$\p^t_H$ rotates each circle at constant speed
\[ \frac{d \p^t_H}{dt} = X_H \circ \p^t_H, \ \ \ \p^0_H = id . \]
Thus $\p^t_H$ is simply the rotation of the 2-sphere about its axis $x_3$ 
through the angle $t$. For any two elements $\b(X_n), X_{n+1}$ in $S^2$, there
is a great circle connecting them with direction $l$ which
is the normal direction to the disk spanned by the great circle. Hence
the rotation about the axis $l$ generates a Hamiltonian flow
${\e}_t$ which connects $\b(X_n), X_{n+1}$ at $t_0$ and $t_1$. Using this
Hamiltonian flow, we get
\[\psi_t(X_1, \cdots, X_{n+1}, X_1, \cdots, X_{n+1})  \]\[ =
(X_1, \cdots, X_{n}, {\e}_t(X_{n+1}), \b(X_1), \cdots, 
\b(X_{n-1}), \b(X_n) {\e}_t(X_{n+1}) \b(X_n)^{-1}, \b(X_n)) \]
from $H_{n+1} \to H_{n+1}$, clearly it satisfies (\ref{product}) and
commutes with $SU(2)$ action. Note that the image of a reducible point
in $H_{n+1}$ under $\psi_t$ is also reducible. So $\psi_t$ maps 
$S_{n+1}$ to itself.
Thus we get a Hamiltonian isotopy
$\psi_t : \hat{H}_{n+1} (= (H_{n+1} \setminus S_{n+1})/SU(2)) \to \hat{H}_{n+1}$
between $\psi_{t_0} = \hat{g}(\p_{\b})$ by (\ref{pb}) and
$\psi_{t_1} = \p_{\s_n \b}$ by (\ref{pnb}). So there is a natural isomorphism
\begin{equation} \label{type2}
HF_i^{\mbox{sym}}(\p_{\s_n \b}) \cong HF_i^{\mbox{sym}}(\hat{g}(\p_{\b})) = 
HF_i^{\mbox{sym}}(\p_{\b}), \ \ \ i \in Z_{2N} .
\end{equation}
The first isomorphism is from the invariance property of the symplectic 
Floer homology under the Hamiltonian isotopy $\psi_t$ and the second
from the natural identification. We can similarly prove that
\[HF_i^{\mbox{sym}}(\p_{\s_n^{-1}\b}) \cong 
HF_i^{\mbox{sym}}(\p_{\b}), \ \ \ i \in Z_{2N} .\]
Combining (\ref{type1}), (\ref{type2}) and the above discussion, we obtain
the desired result. \qed
\begin{cor} \label{sign}
As in Theorem~\ref{invar}, we have the Euler characteristic 
$\chi (HF_*^{\mbox{sym}}(\p_{\b}))$ of
the symplectic Floer homology
\[ \chi (HF_*^{\mbox{sym}}(\p_{\b})) =  - \lc (K) =  \frac{1}{2} sign (K) .\]
\end{cor}
The sign in Corollary~\ref{sign} will be fixed by an example in \S 5. 
For a knot $K$, its Gordian number $u(K)$ (also called the unknotting number)
is the smallest number of times that the string must be allowed to pass
through itself if the knot is to be changed to the unknot.
Murasugi (Theorem 10.1 in \cite{mur}) showed that the absolute value 
of the signature of a knot is not greater than twice the unknotting number. 
Hence combining with Corollary~\ref{sign}, we have
\begin{equation} 
| \chi (HF_*^{\mbox{sym}}(\p_{\b})) | \leq u(K) , \ \ \
\mbox{for $K = \overline{\b}$}.
\end{equation}
The equality holds for the $(2,q)$ torus knot (see \S 5).
We end this section by several remarks which we will discuss them in details
elsewhere.

\noindent{\bf Remarks:} (1) It is natural to ask if there is an instanton Floer
homology of braids which is the original motif of the present paper. In 
particular, based on the work in \cite{he}, one may construct an instanton Floer
homology (like in \cite{ll1}) which is isomorphic to the one we defined
in this section. Note that there is a period 4 instanton Floer homology of knots
in \cite{fl1} and \cite{bd} Part II.

(2) The set-up in \cite{clm} provides a possible relation between the instanton
Floer homology of integral homology 3-spheres and the symplectic Floer homology
of braids.

(3) In \cite{mu}, D. Mullins discussed the Casson-Walker invariant  for
2-fold branched covers of $S^3$. His identification
leads to the relations among the Casson-Walker invariant 
$\lam_{CW}(\hat{S}^3, K)$ of a knot in the rational homology 3-sphere 
$\hat{S}^3$ which is double branched cover of $S^3$ along $K$, 
the signature of the knot $sign(K)$ and the Jones polynomial $V_K(t)$:
\[ 2 \lam_{CW} (\hat{S}^3, K) - \lam_{CL} (K) = \frac{1}{3} \frac{d}{dt}
(\log V_K( -t))|_{t =1} . \]
From our approach, it may shed a light on finding out whether any link between
the Floer homology and the Jone polynomials, one of the most important
problems in the 3-manifold topology.

(4) The computation of $HF_*^{\mbox{sym}}(\p_{\b})$, for $\overline{\b } =K$,
is interesting (and hard) problem in its own right. For symplectic (weakly) 
monotone manifolds and symplectic diffeomorphisms which are exact, 
the Floer homology groups are computed by Floer \cite{fl3}, Hofer and 
Salamon \cite{hs}. For symplectic diffeomorphisms which are isotopic to the
identity through symplectic diffeomorphisms, the Floer homology is computed
by L\^{e} H\^{o}ng and Ono \cite{vo} associated with Calabi invariant. 
It is unknown how to compute the Floer homology
for general symplectic diffeomorphisms, even the Euler characteristic number,
see the Problem 2 in \cite{fl1}. We will develop some techniques to compute the
symplectic Floer homology for the connected sum of knots representing by braids
in a future paper.

\section{Computations of the symplectic Floer homology of braids}

In this section, we mainly concern a computation of
 the symplectic Floer homology
for the trefoil knot in order to fix the sign in Corollary~\ref{sign}.

Let $K_n$ be a subset of $B_n$ which are the representatives of knots.
From the natural inclusion $B_n \hookrightarrow B_{n+1}$, we have a direct 
system of $\{B_n\}_{n \geq 2}$ and a similar direct system
$\{K_n \}_{n \geq 2}$. The direct limit of the direct system 
$\{K_n \}_{n \geq 2}$ exists and
\[ {\cal K} = \lim_{\to} K_n = \cup_{n\geq 2}K_n , \]
is the space of all knots in terms of braids by a theorem of Alexander.
Note that ${\cal K}$ is again a subset (non group structure) of the
direct limit of $\{B_n\}_{n \geq 2}$, ${\cal B} = \lim_{\to }B_n$.
What we have shown in \S 4 is that there is a map $HF$ such that
\begin{equation} \label{new}
HF: {\cal K} \to {\cal HF} ; \ \ \ \ \
\overline{\b} (= K \in K_n) \mapsto 
HF_*^{\mbox{sym}}(\p_{\b}) , \end{equation}
provides a new invariant for knots.

For a unknotted knot $K_0 = \overline{\s_1^{\pm}}$, we have 
\[{\CR}^*(S^3 \setminus K_0)^{[i]} = \emptyset \ \ \ (\mbox{empty set}), 
\ \ \ {\CS}^*(S^3 \setminus K_0)^{[i]} = \{s\} ,\]
by the unknotting theorem. So the Floer chain complex is trivial and
$HF_*^{\mbox{sym}}(\p_{\s_1^{\pm}}) = 0$. By Markov's result, any unknotted
knot can be obtained by finite sequence of Markov moves from $\s_1^{\pm}$, 
then by Theorem~\ref{invar}, 
\begin{equation} \label{unknot}
HF_*^{\mbox{sym}}(\p_{\s}) = 0, \ \ \ \ \
\mbox{ if $\overline{\s}$ is unknotted} .\end{equation}

For the right handed trefoil knot $3_1 = \overline{\s^3_1}$ (or $(2,3)$ torus 
knot), we have that ${\CR}^*(S^2 \setminus \overline{\s^3_1})^{[i]}$
is a 2-sphere with four cone points deleted by \cite{lin} Lemma 2.1. Thus
\[ \pi_1 ({\CR}^*(S^2 \setminus \overline{\s^3_1})^{[i]}) = Z^3, \ \ \
\pi_2({\CR}^*(S^2 \setminus \overline{\s^3_1})^{[i]}) = 0 , \]
the grading for the symplectic Floer homology is an integral grading 
($N \equiv 0$). By \S 4, we use the special reducible {\re}
$\{s\}$ to fix the grading. Recall in \cite{lin}, up to conjugate, we can assume
that
\[ X_1 = \left( \begin{array}{cc}
i \cos \th_1 & \sin \th_1 \\
- \sin \th_1 & - i \cos \th_1  \end{array} \right), \ \ \
X_2 = \left( \begin{array}{cc}
i & 0 \\
0 & -i \end{array} \right), \ \ \ 0 \leq \th_1 \leq \pi . \]
Conditions $X_1 X_2 = Y_1 Y_2$ and $tr(Y_i) = 0, i = 1, 2$ provide
\[ Y_1 = \left( \begin{array}{cc}
i \cos \th_2 & \sin \th_2 \\
- \sin \th_2 & - i \cos \th_2 \end{array} \right),
\ \ \ - \pi \leq \th_2  \leq \pi . \]
So ${\CR}^*(S^2 \setminus \overline{\s^3_1})^{[i]}$ is parameterized by 
$(\th_1, \th_2)$ modulo the involution $(\th_1, \th_2) \to (- \th_1, - \th_2)$.
The orientation is given by the counterclockwise orientation of 
$(\th_1, \th_2)$-plane.
\[{\CR}^*(D^3_- \setminus \overline{\s^3_1})^{[i]} = \{ (\th_1, \th_1) | 
\ \ 0 < \th_1 < \pi \} , \]
\[{\CR}^*(D^3_+\setminus \overline{\s^3_1})^{[i]} = 
Graph(\p_{\overline{\s^3_1}}) = \{ (X_1, X_2, \s^3_1(X_1), \s^3_1(X_2))\},\]
where $Y_1 = \s^3_1(X_1) = (X_1X_2)^2X_1^{-1}(X_1X_2)^{-1}$. A
straightforward calculation shows that
\[ Graph(\p_{\overline{\s^3_1}}) = \{ (\th_1, 4 \th_1) | 
\ \ 0 < \th_1 < \pi \} , \]
\begin{equation}
Fix(\p_{\overline{\s^3_1}}) = \{\r = (\frac{2 \pi}{3}, \frac{2 \pi}{3}) \} = 
{\CR}^*(D^3_+\setminus \overline{\s^3_1})^{[i]} \cap 
{\CR}^*(D^3_- \setminus \overline{\s^3_1})^{[i]}, \end{equation}
is a single element and nondegenerate as in Figure 1.
\medskip

\[ \mbox{Figure 1. The right handed trefoil knot $\overline{\s^3_1}$}\]

The Maslov index of $\r$, $\mu (\r) = \mu (\r, s)$, is the same Maslov index
of two Lagrangian submanifolds
$L_0 = {\CR}^*(D^3_- \setminus \overline{\s^3_1})^{[i]}$ and
$L_1 = Graph(\p_{\overline{\s^3_1}})$ by \cite{clm2}.
This index can be computed as the winding number of a path of $L_1$ about the
diagonal $L_0$ as in \cite{ar} (Note that in \cite{clm2}, $L_0$ is the
$\{(\th_1, 0)\}$). So we have that $L_1$ clockwisely winds the $L_0$ once.
$\mu (\r ) = -1$. Hence
\begin{equation}
HF_i^{\mbox{sym}}(\p_{\s^3_1}) = \left \{
\begin{array}{ll}
Z & \mbox{if $i = -1$} \\
0 & \mbox{if $i \neq -1$} \end{array} \right.
\end{equation} 
The Euler number $\c (HF_*^{\mbox{sym}}(\p_{\overline{\s^3_1}})) = -1 
= - \lc (\overline{\s^3_1})$ in \cite{lin} p 348. Note that the orientation
we used is the opposite one in \cite{ro} p 220. This example fixes
the sign in Corollary~\ref{sign}.

Similarly, we can describe the graph $\p_{\s_1^{-3}}$ for the left
handed trefoil knot $- 3_1 = \overline{\s_1^{-3}}$:
\[ Graph(\p_{\s_1^{-3}}) = \{ (\th_1, \frac{\pi}{2} + 4 \th_1) | \ \ 
0 < \th_1 < \pi \} . \]
The Maslov index of $\r_-$ is zero since there is a $\frac{\pi}{2}$-shift.
\begin{equation}
HF_i^{\mbox{sym}}(\p_{\s_1^{-3}})= \left \{
\begin{array}{ll}
Z & \mbox{if $i = 0$} \\
0 & \mbox{if $i \neq 0$}\end{array} \right. \end{equation}

For the $(2, q)$ torus knot $\overline{\s_1^q}$ ($q$ must be odd), we have
the same method to calculate, by induction,
\[ Graph \p_{\s_1^q} = \{ ( \th_1, (q+1) \th_1 ) | \  \ 0 < \th_1 < \pi \},\]
\[Fix(\p_{\s_1^q}) = \{ \r_k = (\frac{2 \pi k}{q}, \frac{2 \pi k}{q}) | \ \ 
k = 1, 2, \cdots, \frac{q-1}{2} \} . \]
The Maslov indexes are determined by
$\mu(\r_1) = -1$. $\mu(\r_{k-1}, \r_k)$ is the Maslov index of two
Lagrangian submanifolds intersecting transversally at two
smooth points $\r_{k-1}, \r_k$. By \cite{ar} \S 1.4, $\mu(\r_{k-1}, \r_k)$
is the number of rotations of $Det^2$, 
thus $\mu(\r_{k-1}, \r_k) = -2$ (orientation). By additivity of the Maslov
index, we have the nontrivial Floer chain groups:
\[ C_{-2k +1}^{\mbox{sym}}(\p_{\s_1^q}) = Z < \r_{k } >, \ \ \
 k = 1, 2, \cdots, \frac{q-1}{2} .\]
Therefore we have the trivial Floer boundary for this case again. So
\begin{equation}
HF_i^{\mbox{sym}}(\p_{\s_1^q})= \left \{
\begin{array}{ll}
Z & i = - 2k +1, \ \ k = 1, 2, \cdots, \frac{q-1}{2}\\
0 & \mbox{otherwise}
\end{array} \right. \end{equation}
Hence $\chi (HF_*^{\mbox{sym}}(\p_{\s_1^q})) = - \frac{q-1}{2}$,
the signature of the $(2,q)$ torus knot is $- (q-1)$ in \cite{lit}.
The number $\frac{q-1}{2}$ is also the unknotting number of 
the $(2,q)$ torus knot by \cite{mur}.

In order to compute the Floer homology, 
one has to study the possible nontrivial boundary map. This is the case
for general $(p, q)$ torus knot (a braid representative $(\s_1 \s_2 \cdots
\s_{p-1})^q, (p,q) = 1$). There are some methods in determining
the {\re}s of knot groups in $SU(2)$ (not using braid {\re}s) as in
\cite{bu} and \cite{kl}.

\noindent{\bf Acknowledgment:} We presented this paper for the first
time at 911TH AMS meeting at Baton Rouge. I would like to thank
the organizers P. Gilmer, R. Litherland and N. Stoltzfus for the
invitation. I am very grateful to R. Lee for many interesting and helpful
conversations and suggestions on this work. I would also like to thank 
my colleagues W. Jaco and R. Myers for many discussions.

\end{document}